\documentclass{arxiv}
\pdfoutput=1

\usepackage{graphicx}
\usepackage{mathrsfs}
\usepackage{verbatim}

\usepackage[utf8]{inputenc}
\mleftright
\usepackage{subcaption}

\newcommand{\norm}[1]{\left\lVert#1\right\rVert}

\usepackage{subcaption}
\begin{document}

\title[Signatures in Shape Analysis]{Signatures in Shape Analysis: an Efficient Approach to Motion Identification}
%
%
\author{ Elena Celledoni}
\email{elena.celledoni@ntnu.no}
\urladdr{https://www.ntnu.edu/employees/elena.celledoni}
\author{Pål Erik Lystad}
\email{paalel@stud.ntnu.no}
\address{Department of Mathematical Sciences, Norwegian University of Science and Technology, 7491 Trondheim, Norway.}
\author{Nikolas Tapia}
\email{tapia@wias-berlin.de}
\urladdr{https://www.wias-berlin.de/people/tapia}
\address{Weierstraß-Institut Berlin, Mohrenstr. 39, 10117 Berlin, Germany.}
%
%
%
\maketitle              

\begin{abstract}
Signatures provide a succinct description of certain features of paths in a reparametrization invariant way. We propose a method for classifying shapes based on signatures, and compare it to current approaches based on the SRV transform and dynamic programming.

\smallskip
\noindent
{\textbf{\keywordsname.} Shape analysis, Signature, Motion identification}
\end{abstract}




\section{Introduction}
    Shape analysis is a broad and growing subject addressing the analysis of different types of data ranging from  surfaces, landmarks, animation data etc. 
    In this paper shapes are unparametrized curves. Mathematically a shape is an equivalence class of curves under reparameterization, that is, two curves $c_0, c_1:[0,1]\to M$ are equivalent and  determine the same shape if there exists a strictly increasing smooth bijection $\varphi:[0,1]\to[0,1]$ such that $c_1 = c_0\circ \varphi$.
    For a given curve $c$ we denote by $[c]$ the corresponding shape.

    The similarity between two shapes $[c_0],[c_1]$ is then defined by creating a distance function $d_\mathcal{S}$ on the space of shapes $\mathcal{S}$,
    \begin{equation}\label{eq:intro-distance-measure}
        d_\mathcal{S}([c_0],[c_1]) := \inf_{\substack{\varphi}}d_\mathcal{P}(c_0, c_1\circ\varphi)
    \end{equation}
    where $d_{\mathcal P}$ is a suitable reparameterization invariant Riemannian distance on the manifold of parametrized curves.

    Finding the optimal reparameterization $\varphi$ is however computationally demanding, and in many applications simply unnecessary. This is specifically the case of applications where the optimal parametrization is not explicitly used for further calculations, e.g. problems of identification and classification. Ways of circumventing this step are therefore of great interest.
    
    In recent years, after extensive work by Terry Lyons and collaborators, the theory of rough paths has gained considerable importance as a toolbox for mathematical analysis and for mathematical modeling in applications. In this context, the signature map provides a faithful representation of paths, capturing their essential global properties. A fundamental property of the signature is its invariance under reparameterization, surmising its importance for shapes.
    
    In this paper, we define a measure of similarity between shapes in $\mathcal{S}$ by means of the signature. We define a distance directly on $\mathcal{S}.$
    We test the viability of this approach and use it to classify motion capture animations from the CMU motion capture database \cite{CMUGraphicsLabMotionCaptureDatabase}. Indeed, this leads to an efficient technique that delivers results comparable to what is obtainable with methodologies based on the SRV transform, but at a much lower computational cost.

\section{Shape analysis on Lie groups}\label{section:shape-analysis-framwork}
    In the following, $G$ will denote a finite-dimensional Lie group under multiplication with identity element denoted by $e$.
    We let $\mathfrak{g}$ denote the corresponding right Lie algebra $\mathfrak{g}:=\mathcal{L}_R(G)$. 
    For a fixed $g\in G$, left and right translation by $g$ will be denoted $L_g(h) = g\cdot h$ and $R_g(h) = h\cdot g$ respectively.
    
    \subsection{Shape Space}
     We consider the space $C^\infty ([0,1],G)$ of parameterized smooth curves on $G$, i.e. smooth maps $c:[0,1]\to G$.
     To model the curves as unparameterized, or independent of parameterization, we define the \emph{shape space} $\mathcal{S}$ as the quotient space
    \begin{equation}\label{eq:shape-space}
    \mathcal{S} = C^\infty([0,1],G) \slash\text{Diff}^+,
    \end{equation}
    where $\text{Diff}^+$ is the group of orientation preserving diffeomorphisms of the parameter space $[0,1]$.
    The elements of $\mathcal{S}$ are equivalence classes of curves. The elements of the same class are curves which can be mapped to one another by changing their parameterization, that is, two curves $c_0,c_1 \in C^\infty (I, G)$ are equal in shape space if there exists $\varphi\in\operatorname{Diff}^+$ such that $c_1=c_0\circ\varphi$.
     
    In the setting of our application, the search for optimal time parametrizations can be viewed as syncing up the animations, removing disturbances due to small pauses, different periodicity, or asynchronous starting and stopping, by shifting the movement of one character to match the other as closely as possible.

\subsection{Geodesic Distances on Shape Space}\label{section:geodesic-distance-on-shape-space}

    Our goal is to introduce a meaningful and computable distance $d_\mathcal{S}$ on $\mathcal{S}$ to estimate the similarity between two shapes. This area of research started with the efforts of Younes \cite{Younes98computableelastic}.
    We will restrict the space of curves to the space of immersions, i.e. curves with non-vanishing first derivative, which we denote by
    \begin{equation}\label{eq:immersion-space}
        \mathcal{P} = \text{Imm}([0,1], G).
    \end{equation}
  
  
    Let $d_\mathcal{P}$ be a pseudo-metric on $\mathcal{P}$. We define $d_\mathcal{S}$, for two elements $[c_0],[c_1]\in \mathcal{S}$, by
    \begin{equation}\label{eq:metric-shape_space}
        d_\mathcal{S}([c_0],[c_1]) := \inf_{\varphi \in \operatorname{Diff}^+} d_\mathcal{P}(c_0, c_1 \circ \varphi).
    \end{equation}
    As shown in \cite[Lemma 3.4]{CelledoniShapeLieComputerAnimation}, $d_S$ will be a pseudo-metric on $\mathcal{S}$ if $d_\mathcal{P}$ is a \textit{reparameterization invariant} or, in other words, if for any two $c_0,c_1 \in \mathcal{P}$ and any $\varphi\in\operatorname{Diff}^+$ we have that
    \begin{equation}\label{pseudo-matric-reparam-invariant}
        d_\mathcal{P}(c_0 \circ \varphi, c_1 \circ \varphi) = d_\mathcal{P}(c_0,c_1).
    \end{equation}
    
    An obvious choice of metric on $\mathcal P$ is the familiar $L_2$-metric. However, as shown by Michor and Mumford \cite{MichorMumfordVanishingGeodesic}, this metric leads to vanishing geodesic distance which renders it useless. 
    They further show in \cite{MichorMumforRiemannianMetricHamiltonian} that one solution to this problem is to consider metrics based on arc-length derivatives, creating a class of Sobolev-type metrics.
    
    There are multiple possible metrics in this class. One option is based on what is usually referred to as the \emph{Square Root Velocity Transform} (SRVT). This transform and accompanying metric was first introduced, in the context of shape analysis, by Srivastava et al. \cite{SrivastavaKlassenShapeAnalysisElasticCurves}, who used the transformation when working with curves in Euclidian spaces. The transformation has later been adopted to more general shapes. 
    Of particular interest is the formulation for shapes that are represented as Lie-group valued curves \cite{CelledoniShapeLieComputerAnimation}.
    
    We define the SRVT $\mathcal{R}:\mathcal P\to C^\infty([0,1],\mathfrak g\setminus\{0\})$ by
    \begin{equation}\label{eq:SRVT}
            \mathcal{R}(c)(t)\coloneq\frac{R^{-1}_{c(t)*}(\Dot{c}(t))}{\sqrt{\norm{\Dot{c}(t)}}}.
    \end{equation}
    This transformation has the following useful properties \cite[Lemma 3.6]{CelledoniShapeLieComputerAnimation}:
    \begin{enumerate}
        \item For every $c\in\mathcal{P}$ and $\varphi\in\operatorname{Diff}^+$, the following equivariant property holds:
        \begin{equation}\label{eq:reparameterization-relation-srvt}
            \mathcal{R}(c\circ\varphi) = \mathcal{R}(c)\circ\varphi\cdot\sqrt{\dot\varphi}.
        \end{equation}
        \item It is translation invariant: for all $c\in\mathcal{P}$ and $g\in G$
        \begin{equation*}
            \mathcal{R}(R_g(c)) = \mathcal{R}(c).
        \end{equation*}
    \end{enumerate}
    A similar result is true for shapes with values in Euclidean spaces \cite{SrivastavaKlassenShapeAnalysisElasticCurves}.

    
    Further, one can obtain a Riemannian metric $d_{\mathcal{P}*}$ that coincides with the geodesic distance on a submanifold $\mathcal{P}_* \subset\mathcal{P}$ by using the SRVT to pull back the $L_2$-metric on $C^\infty(I,\mathfrak{g}\setminus\{0\})$ \cite{CelledoniShapeLieComputerAnimation}.
    Further restricting the immersion space to $\mathcal{P}_* = \{c \in \mathcal{P}:c(0)=e\}$, where $e$ is the identity element in $G$, the distance $d_{\mathcal{P}_*}$ turns out to be reparameterization invariant.
    
    This invariance implies, in particular, that it will also yield a geodesic distance on $\mathcal{S}_*:=\mathcal{P}_*/\text{Diff}^+$ \cite{BruverisGeodesicCompletenessSobolev}. The restriction to $\mathcal{P}_*$ isn't very troublesome as any curve can be transferred to this space by right translation by the inverse of its initial value, that is $R_{c(0)^{-1}}$ \cite{CelledoniShapeLieComputerAnimation}. 

    Using the equivariant property for the SRVT from equation \eqref{eq:reparameterization-relation-srvt} and defining $q_i=\mathcal{R}(c_i)$ for $i=0,1$, the problem of calculating the metric for the shape space $\mathcal S_*$ in equation \eqref{eq:metric-shape_space} can be written as
    \begin{equation}\label{eq:optimization-problem-metric-formulation}
        d_{\mathcal{S}_*}(c_0,c_1) = \inf_{\substack{\varphi \in \text{Diff}^+(I)}} \sqrt{\int_I{\Vert q_0(t)-q_1(\varphi(t))\cdot\sqrt{\dot{\varphi}}\Vert^2 dt}}.
    \end{equation} 
    
    Finding this infimum will generally be very difficult. The usual approach is therefore to discretize the curves and solve instead a finite  dimensional optimization problem. The most common methods used to solve this problem in shape analysis \cite{SrivastavaKlassenShapeAnalysisElasticCurves} are based on either the gradient descent method or a dynamic programming algorithm (DP). In our experiments we  use the DP approach described in \cite{BauerLandmark-GuidedElasticShapeAnalysis}.

\section{Signatures}

Signatures, introduced by K.-T. Chen \cite{KTChenSig} for smooth paths and later generalized by Lyons \cite{Lyons1998} under the name of geometric rough paths, are an important tool for the study of the solutions of controlled differential equations, but have also proved useful for solving classification problems of time series, Machine Learning and Topological Data Analysis \cite{ChevyTopoData}.

In the usual framework, signatures are defined for paths taking values in a Banach space.
From a geometric point of view, and in light of our purposes, this setting has to be adapted.
Luckily, Chen also considered signatures for curves taking values on a smooth manifold \cite{KTChenSig}.
This definition is quite general and relies on the selection of a frame bundle.
For Lie groups there is a canonical choice: the Maurer--Cartan form. This is the unique right-invariant one form $\omega$ such that $\omega_e=\mathrm{id}_{\mathfrak{g}}$, i.e. $\omega(v)=(R_g^{-1})_*v$ for $v\in T_g(G)$ \cite[p. 311]{HilgertNeebStructureGeometryLieGrups}.

Below we denote, for a finite-dimensional vector space $V$ of dimension $d=\dim V$, the tensor algebra over $V$,
\[ T(V)\coloneq\bigoplus_{n\ge0}V^{\otimes n}. \]
We observe that $T(V)$ is always infinite-dimensional.
Its dual space is denoted by $T(\mkern-3mu(V)\mkern-3mu)\coloneq T(V)^*$, and it may be identified with the ring of formal power series in $d$ noncommuting variables $\{e_1,\dotsc,e_d\}$.

\begin{dfn}
    Let $G$ be a $d$-dimensional Lie group and $\alpha\in C^\infty([0,1],G)$ be a smooth curve and $\omega$ the Maurer-Cartan form on $G$.
    The signature $S(\alpha)$ of $\alpha$ is the family of linear maps on $T(\mathbb R^d)$ recursively defined by $\langle S(\alpha)_{s,t},1\rangle\coloneq 1$ and
    \[ \langle S(\alpha)_{s,t},e_{i_1\dotsm i_p}\rangle\coloneq\int_s^t\langle S(\alpha)_{s,u},e_{i_1\dotsm i_{p-1}}\rangle\,\omega^{i_p}_{\alpha(u)}(\dot\alpha(u))\,\mathrm du. \]
\end{dfn}
In this definition, the notation $\omega^j_g(v)$ denotes the $j$-th component of the vector $\omega_g(v)\in\mathfrak g$ in a basis of the Lie algebra $\mathfrak g$ of $G$.

The signature provides a compact description of certain features of a path \cite{KTChenIntegrationOfPaths}.
One of its main advantages in our context is its reparameterization invariance: for any orientation-preserving diffeomorphism $\varphi$ on $[s,t]$ we have that $$S(\alpha\circ\varphi)_{s,t} = S(\alpha)_{s,t}.$$

Other fundamental properties include:
\begin{enumerate}
    \item For each $0\le s<t\le 1$, the signature $S(x)_{s,t}$ belongs to the set of group-like elements of $T(\mkern-3mu(\mathbb R^d)\mkern-3mu)$, and for any $0\le s\le 1$, $S(x)_{s,s}=1$, the neutral element in the group.
    \item \textbf{Chen's rule:} For any three $0\le s<u<t\le 1$ we have \[ S(x)_{s,u}\otimes S(x)_{u,t} = S(x)_{s,t}. \]
\end{enumerate}
Using these properties, signatures may be efficiently computed for some restricted classes of paths. For example, if $x$ is a straight line in $\mathbb R^d$ with base point $a\in\mathbb R^d$ direction $b\in\mathbb R^d$, i.e. $x_t=a+t b$ for $t\in[0,1]$, then
\begin{equation}
\begin{split}
    S(x)_{s,t} &= \exp_{\otimes}((t-s)b)\\
    &= 1 + (t-s)b + \frac{(t-s)^2}{2}b\otimes b + \frac{(t-s)^3}{6}b\otimes b\otimes b + \dotsb.
\end{split}
\label{eq:sigexp}
\end{equation}
A similar statement is true for geodesic curves on a finite-dimensional compact Lie group.

We may think of signatures as an infinite vector indexed by \emph{words} over the alphabet $\{1,\dotsc,d\}$. In particular, for a piecewise linear path the above formula means that if we want to know the component in \eqref{eq:sigexp} corresponding to the word $w=i_1\dotsm i_k$ then
\[ \langle S(x)_{s,t},e_w\rangle= \frac{(t-s)^k}{k!}\prod_{j=1}^kb_{i_j}
\]

For a general piecewise linear path $x$, we may use the above formula and Chen's rule to deduce that
\[ S(x)_{s,t} = \exp_\otimes(\Delta t_1 b_1)\otimes\exp_\otimes(\Delta t_2 b_2)\otimes\dotsm\otimes\exp_\otimes(\Delta t_m b_m) \]
where $\Delta t_k=t_k-t_{k-1}$ are the length of the time intervals where the path is sampled and $b_1,\dotsc,b_k$ are the slopes of the path in each of these intervals.
The entries of this expression may be computed by using a Baker--Campbell--Hausdorff-type formula, for example.

Finally, we remark that the signature possesses another interesting property, namely it is an homomorphism from path space with concatenation to the tensor algebra $T(\mkern-3mu(\mathbb R^d)\mkern-3mu)$. This means that if we are given two paths $x\colon[0,1]\to G$ and $y\colon[0,1]\to G$, and we concatenate them to form a new path $x\cdot y$, then
\[ S(x\cdot y)_{0,1} = S(x)_{0,1}\otimes S(y)_{0,1}. \]
Moreover, if we reverse the path $x$, i.e. we define $\overleftarrow{x}(t)\coloneq x(1-t)$ then
\[ S(\overleftarrow x)_{0,1} = S(x)_{0,1}^{-1} \]
where the inverse is taken in the group-like elements of the tensor algebra.

It can be shown that actually, as a function of time the signature satisfies the differential equation
\[ \frac{\mathrm d}{\mathrm dt}S(x)_{s,t}=S(x)_{s,t}\otimes\dot x_t, \quad S(x)_{s,s}=1 \]
in the tensor algebra.
From this point of view, the signature map corresponds to the flow map of the vector field given by the base path.
Thus, the signature belongs to an infinite-dimensional Lie group whose Lie algebra is the free Lie algebra over $\mathbb R^d$ which we denote by $\mathfrak L(\mathbb R^d)$. It does not, however, constitute a one-parameter subgroup.
Therefore, for each fixed time interval $[s,t]$ we can map the signature to the free Lie algebra via a logarithm map, and we define
\[ \Lambda(x)_{s,t}=\log(S(x)_{s,t})\in\mathfrak L(\mathbb R^d). \]
This element, called the log-signature in the literature, provides a minimal description of the path, which is equivalent to the full signature.

There are many ways in which signatures can be used to compare shapes, but the essential feature is that since the map $S$ is reparameterization invariant, one obtains a way of directly comparing shapes instead of parameterized curves.
For our experiments we chose a particular distance on $T(\mkern-3mu(\mathbb R^d)\mkern-3mu)$ (see next section for the precise formula), but this is by no means the only possible choice.

In making this choice one has to truncate the signature to obtain a finite-dimensional object.
Due to the factorial decay of iterated integrals little information is lost in the process; still, some level has to be chosen and usually this done by running experiments.
Once the truncation level is chosen, several choices of metric are available: the truncated tensor algebra becomes finite-dimensional so it has a nice linear structure and we are free to choose norms on it subject to some compatibility restrictions.
There is also the notion of homogeneous norm on group-like elements, which takes into account the geometry of this group.
Finally, the logarithm in this group maps signatures into a linear space (the free Lie algebra) in a bijective way, so no information is lost, but there is a substantial dimensional reduction.

According to our observations, is the last option which represents the most robust choice in terms of noise sensitivity, while also providing an accurate way of comparing signatures.


\section{Experiments}\label{section:experiment}
    Motion capture animations are usually recorded as the angle of every joint in a skeleton for every frame in an animation. A natural setting for the rotating joints is the Lie group of 3D rotations, $SO(3)$. Every frame consists of $d$ independently rotating joints so the frame can be modeled as an element in $SO(3)^d$, where $SO(3)^d$ is the Cartesian product of $d$ copies of $SO(3)$. Interpolating between the frames will then allow us to model the animation as a parameterized curve.
    
    We use an interpolation scheme in which one uses the log map to linearly interpolate on the Lie algebra, and then pull back to the Lie group with the exponential map. 
    Let $A,B\in SO(3)$, we define the interpolation $\kappa:[0,1]\to SO(3)$ between $A$ and $B$ as
    \begin{equation*}
        \kappa(s) \coloneq \exp\left(s\log\left(B \cdot A^T\right)\right) \cdot A.
    \end{equation*}
    Notice that $\kappa(0)=A$ and $\kappa(1)=B$. Applying this interpolation component-wise to the frames in $SO(3)^d$ will enable us to construct a piece-wise interpolation between the frames of the animation. The Maurer--Cartan form along the interpolation is piece-wise constant, making it easy to compute SRV representations, $d_{\mathcal{P}_*}$-metrics, and signatures.

     To test the effectiveness of the proposed frameworks we check whether they are able to identify different types of character motion. We have selected animations from the CMU motion capture database with descriptions "walk", "run/jog" and "forward jump". These are similar in length, and should produce results that conform with human intuition. 
     
     The test will calculate a distance matrix using the proposed similarity measures. From the distance matrix we produce a multidimensional scaling plot(MDS), depicting how similar, or dissimilar, the animations are. MDS tries to place the data points in 2-dimensional scatter plot while preserving the distances given by the distance matrix. See Kruskal \cite{KruskalWishMultiScale1978} for more more information on this method.
     
\begin{figure}[ht!]
       \centering
       \includegraphics[width=0.8\linewidth]{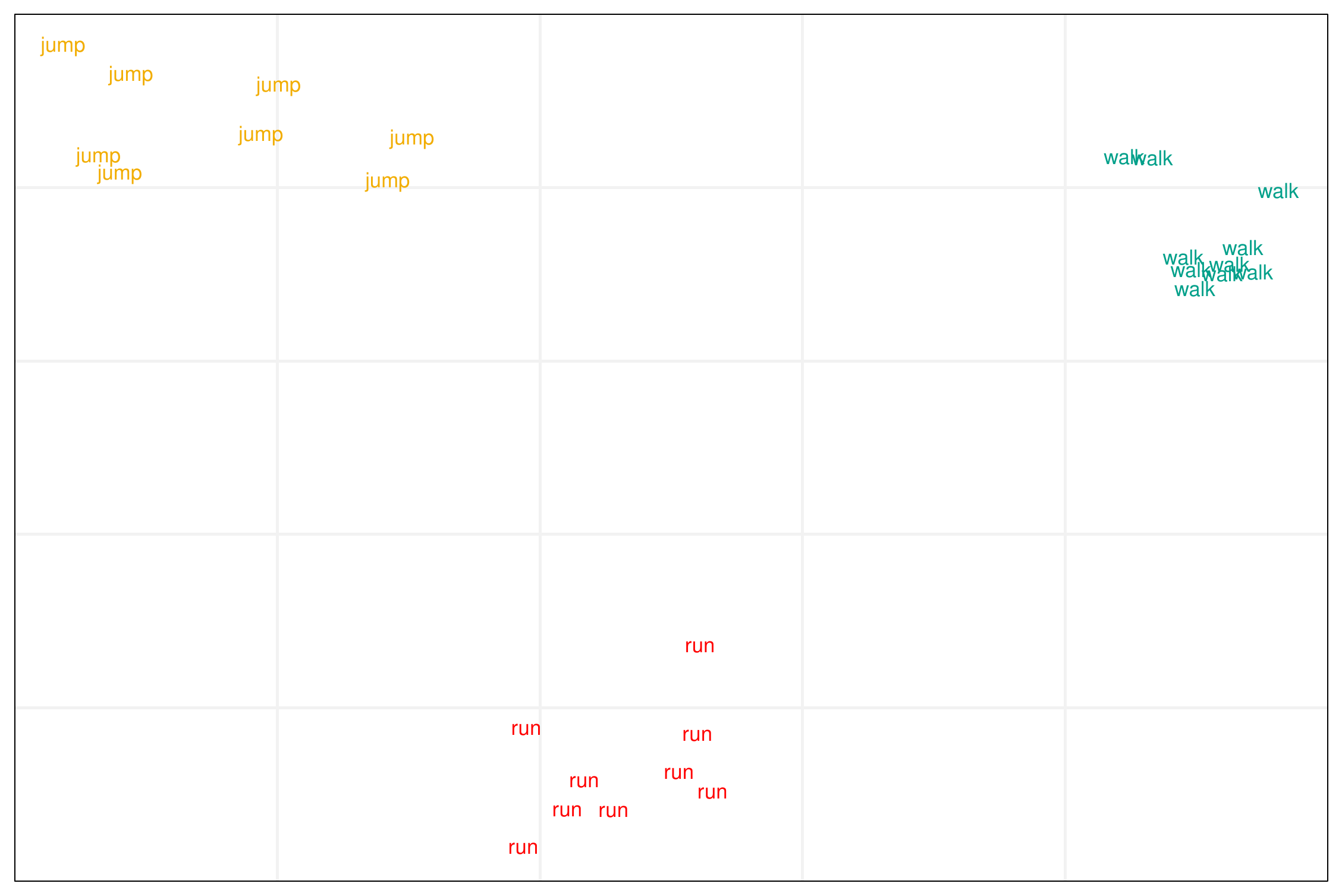}
          \caption{Animations projected to $\mathcal{P}_*$ with distance matrix calculated with the metric $d_{\mathcal{P}_*}$.}
          \label{fig:MDS-regular-run-walk-jump}
      \end{figure}

\begin{figure}[ht!]
        \includegraphics[width=0.8\linewidth]{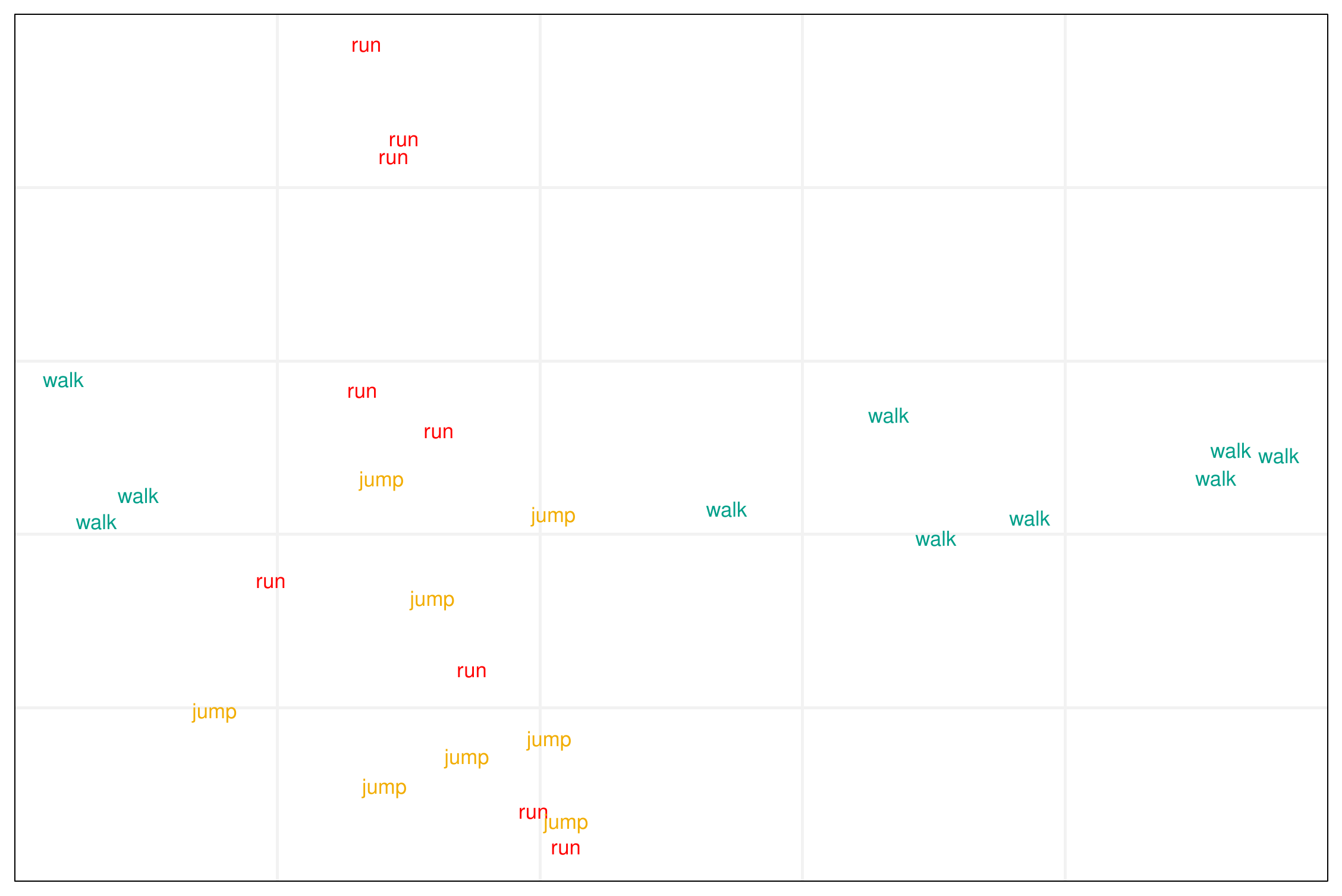}
          \caption{Animations projected to $\mathcal{S}_*$ with distance matrix calculated with metric $d_{\mathcal{S}_*}$ using a DP algorithm.}
          \label{fig:MDS-invariant-run-walk-jump}
    \end{figure}

     In Figure \ref{fig:MDS-regular-run-walk-jump} we calculate the distance matrix using the metric $d_{\mathcal{P}_*}$ on interpolation curves in $\mathcal{P}_*$, and in Figure \ref{fig:MDS-invariant-run-walk-jump} we use the metric $d_{\mathcal{S}_*}$, equation \eqref{eq:optimization-problem-metric-formulation}, on the shapes generated by the curves in $\mathcal{S}_\ast$, where the optimal reparameterization is calculated with a DP algorithm. There are little to no patterns when projecting to the space $(\mathcal{P}_*, d_{\mathcal{P}_*})$, as seen in Figure \ref{fig:MDS-regular-run-walk-jump}. In Figure \ref{fig:MDS-invariant-run-walk-jump} however, we observe that modelling the curves as being parameterization invariant yields three easily distinguishable clusters of animations. Compared to Figure \ref{fig:MDS-regular-run-walk-jump} we see a big benefit from this model assumption.
     
\begin{figure}[h!t]
        \centering
        \includegraphics[width=0.8\textwidth]{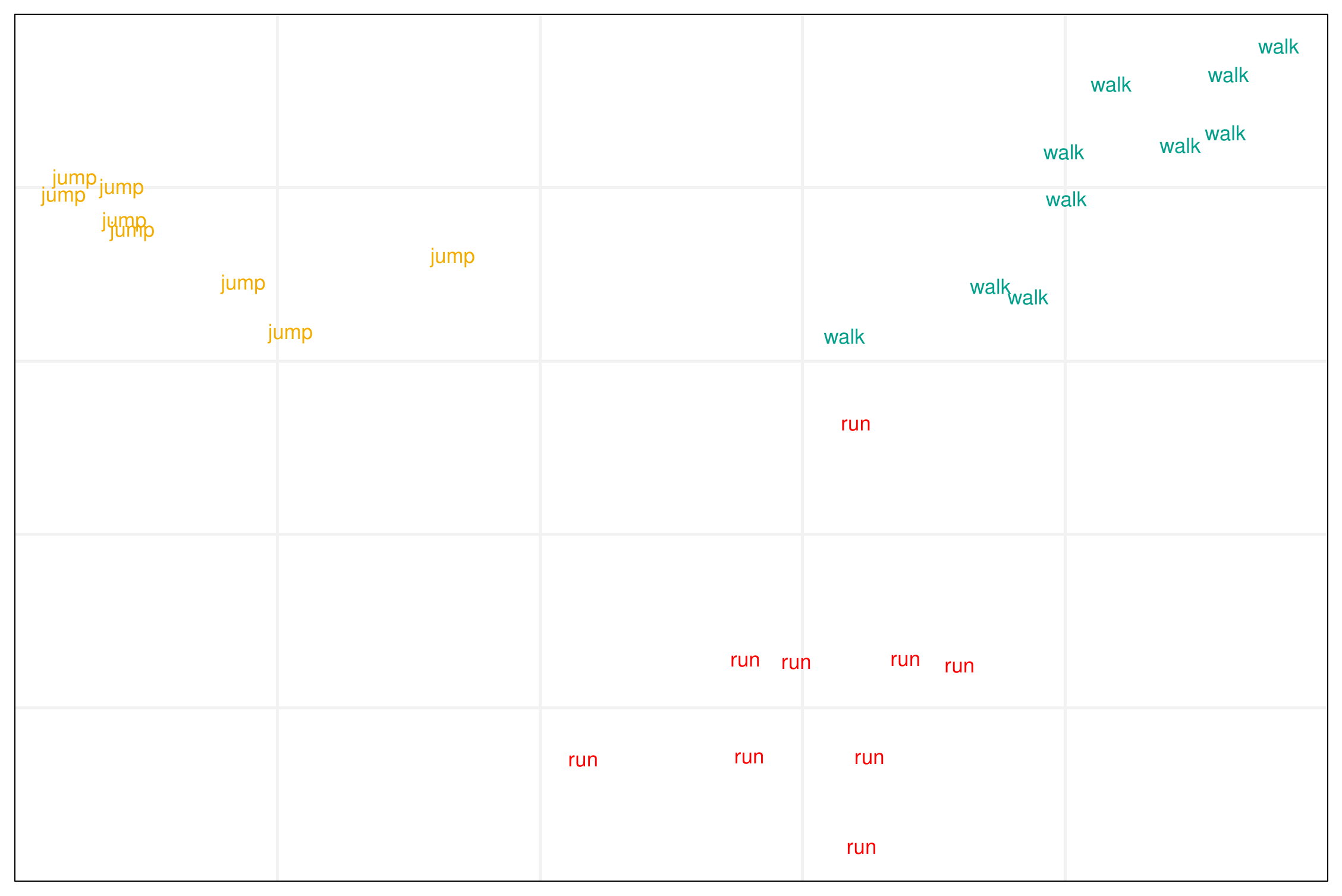}
        \caption{Multi dimensional scaling plot of distance matrix calculated from by projecting animations to the space $\mathcal{S}_\ast$ equipped with the distance function $d_\text{sig}$. In this plot we have taken animation with  descriptions "run/jog", "forward jump" and "walk" from the CMU Motion Capture Database \cite{CMUGraphicsLabMotionCaptureDatabase}.} 
        \label{fig:MDS-signature}
    \end{figure}

    In Figure \ref{fig:MDS-signature} the animations are projected to the shape space $\mathcal{S}$ equipped with the distance function $d_\text{sig}(c_0, c_1) = \norm{\frac{\log S(c_0)}{\norm{\log S(c_0)}}-\frac{\log S(c_1)}{\norm{\log S(c_1)}}}$. While this figure does reveal the same structure as seen in figure \ref{fig:MDS-invariant-run-walk-jump}, the clusters exhibit both a higher internal and a lower external variability. An important take away from this experiment is that this distance function in fact does preserve some of the structure of the shape space.
   

\newpage
\section{Concluding Remarks}
    Our preliminary experiments, show that classifying animations using a distance function on $\mathcal{S}_\ast$ based on signatures produces very encouraging results.
    The proposed method is  computationally very efficient, even though somewhat less accurate than known methods in shape analysis.
  
    The Riemannian metric \eqref{eq:metric-shape_space} requires calculating the optimal reparameterizations between every pair of animations. The proposed signature method instead only requires calculating the signature once for every animation, and then compares animations by computing inexpensive norms. The optimisation procedure is no longer necessary.
    \footnote{A more thorough analysis of the run time complexities associated with these algorithms will be performed in future work.}
    
    In our experiments, the signature method outpreformed the optimal reparameterization metric by a factor of $\sim 2000$ when classifying animations.
    A more precise comparison with the SRVT approach and other methods, see e.g. \cite{LahiriMachinignPLCurves} goes beyond the scope of this work and will be considered in future work. Still our preliminary experiments give an idea of the possible performance benefits gained with the signature approach.
    
    Increasing the accuracy of the signature method might also be possible by defining a more precise similarity measure.
    Nonetheless, our results can be seen as proof of concept for using signatures as an efficient way of classifying shapes.

    \subsection{Acknowledgements}
    This paper contains work done as part of P.E.L.'s master thesis. The master thesis will be published separately as part of NTNU's Master of Science program in Applied Physics and Mathematics \cite{LystadMasterThesis}.
    N.T. acknowledges that part of this work was carried out during his tenure of an ERCIM `Alain Bensoussan' Fellowship Programme (contract number 2018-10) at NTNU. This work was supported by the European Union's Horizon 2020 research and innovation programme under the Marie Sklodowska-Curie, grant agreement No.691070. 
    
    The data used in this project was obtained from \url{http://mocap.cs.cmu.edu}. The database was created with funding from NSF EIA-0196217.

\bibliographystyle{arxiv}
\bibliography{bibliography}


\end{document}